\documentclass[12pt]{article}

\usepackage{amsfonts}

\newtheorem{thm}{Theorem}[section]
\newtheorem{lemma}[thm]{Lemma}

\newtheorem{definition}[thm]{Definition}
\newtheorem{prop}[thm]{Proposition}

\newcommand{\proof
}{\par\medskip\noindent {\bf Proof.\ \ }}

\newcommand{\be}{\begin{equation}}
\newcommand{\ee}{\end{equation}}
\newcommand{\openbox}{\leavevmode
  \hbox to8pt{\hfil\vrule\vbox to6pt{\hrule width6pt\vfil\hrule}\vrule}}

\newcommand{\qed}{\hbox to5pt{ } \hfill \openbox\bigskip\medskip}

\newcommand{\Fq}{\mathbb F _q}

\newcommand{\cS}{\mbox{$\cal S$}}

\newcommand{\cF}{\mbox{$\cal F$}}

\newcommand{\N}{\mathbb N}

\newcommand{\F}{\mathbb F}

\title{Inequalities for two systems of subspaces with prescribed intersections}
\author{G\'abor Heged\"{u}s
\\{\normalsize Antal Bejczy Center For Intelligent Robotics} 
\\{\normalsize Kiscelli utca 82, Budapest, Hungary, H-1032}
\\{\normalsize hegedus.gabor@nik.uni-obuda.hu}
}

\begin{document}



\maketitle

\begin{abstract}
Let $W$ denote a linear space over a fixed field ${\F}$. We define the notions of  weak $ISP$-system and weak $(u,v)$-system $\cS=\{(U_i,V_i):~ 1\leq i\leq m\}$ of subspaces of $W$. We give upper bounds for the size of  weak $ISP$-systems and weak $(u,v)$-systems.
\end{abstract}

\section{Introduction}

First we recall the notion of $q$-binomial coefficients.

The {\em $q$-binomial coefficient} $\Big[ { n \atop m} \Big]_q$ is a $q$-analog for the binomial coefficient, also called a Gaussian coefficient or a Gaussian polynomial. The  $q$-binomial coefficient is given by 
\begin{equation} \label{Gaussco}
 \Big[ { n \atop m}\Big]_q:=\frac{[n]_q!}{[n-m]_q!\cdot[m]_q!}
\end{equation}
for $n,m \in \N$, where $[n]_q!$ is the $q$-factorial (see \cite{K}, p. 26)
$$
[n]_q!:=(1+q)\cdot (1+q+q^2)\cdots (1+q+q^2+\ldots +q^{n-1}).
$$
Clearly we have $\Big[ { n \atop k} \Big]_q=\Big[ { n \atop n-k} \Big]_q$. If we substitute $q=1$ into (\ref{Gaussco}), then this substitution reduces this definition to that of binomial coefficients.
 
Bollob\'as proved in \cite{B} the following two remarkable results in extremal combinatorics.
\begin{thm} \label{Boll}
Let $A_1, \ldots A_m$ and  $B_1, \ldots B_m$ be finite sets satisfying the conditions 
\begin{itemize}
\item[(i)] $A_i \cap B_i =\emptyset$ for each $1\leq i \leq m$;
\item[(ii)] $A_i\cap B_j\ne \emptyset$ for each $i\ne j$ ($1\leq i, j \leq m$).
\end{itemize}
Then 
$$
\sum_{i=1}^m \frac{1}{{|A_i|+|B_i| \choose |A_i|}}\leq 1.
$$
\end{thm}

\begin{thm} \label{Boll2}
Let $A_1, \ldots A_m$ be $r$-element sets and  $B_1, \ldots B_m$ be $s$-element sets such that 
\begin{itemize}
\item[(i)] $A_i \cap B_i =\emptyset$ for each $1\leq i \leq m$;
\item[(ii)] $A_i\cap B_j\ne \emptyset$ for each $i\ne j$ ($1\leq i, j \leq m$).
\end{itemize}
Then
$$
m\leq {r+s \choose s}.
$$
\end{thm}

Tuza proved the following two versions of Bollob\'as Theorem.
\begin{thm} \label{Tuza1}
Let $p$ be an arbitrary real number, $0<p<1$ and $t:=1-p$. 

Let $A_1, \ldots A_m$ and  $B_1, \ldots B_m$ be finite sets satisfying the conditions 
\begin{itemize}
\item[(i)] $A_i \cap B_i =\emptyset$ for each $1\leq i \leq m$;
\item[(ii)] $A_i\cap B_j\ne \emptyset$ or $A_j\cap B_i\ne  \emptyset$ for $i\ne j$ ($1\leq i, j \leq m$).
\end{itemize}
Then
$$
\sum_{i=1}^m p^{|A_i|}t^{|B_i|}\leq 1.
$$
\end{thm}
\begin{thm} \label{Tuza2}
Let $A_1, \ldots A_m$ be $r$-element sets and  $B_1, \ldots B_m$ be $s$-element sets satisfying the conditions 
\begin{itemize}
\item[(i)] $A_i \cap B_i =\emptyset$ for each $1\leq i \leq m$;
\item[(ii)] $A_i\cap B_j\ne \emptyset$ or $A_j\cap B_i\ne  \emptyset$ for $i\ne j$ ($1\leq i, j \leq m$).
\end{itemize}
Then
$$
m\leq \frac{(r+s)^{r+s}}{r^r s^s}.
$$
\end{thm}

Z. Tuza raised in \cite{T} the following question: Let $a,b$ be fixed positive integers. Determine the largest integer $m:=m(a,b)$ such that 
there exists a system  $\cS=\{(A_i,B_i):~ 1\leq i\leq m\}$  of $m(a,b)$ pairs of sets satisfying the  conditions:
\begin{itemize}
\item[(i)] $A_1, \ldots A_m$ are $r$-element sets and  $B_1, \ldots B_m$ are $s$-element sets;
\item[(ii)] $A_i \cap B_i =\emptyset$ for each $1\leq i \leq m$;
\item[(iii)] $A_i\cap B_j\ne \emptyset$ or $A_j\cap B_i\ne  \emptyset$ for $i\ne j$ ($1\leq i, j \leq m$).
\end{itemize}
Tuza proved the following properties of the numbers $m(a,b)$ in \cite{T}.

\begin{prop} \label{Tuza15}
$m(a,1)=2a+1$ for each $a\geq 1$. For every $a,b,\geq 1$
$$
m(a,b)\geq m(a,b-1)+m(a-1,b).
$$
\end{prop}

Proposition \ref{Tuza15} gives a lower bound for $m(a,b)$ near to $2{a+b \choose a}$ for every $a$ and $b$.

Lov\'asz used in \cite{L1} tensor product methods to prove the following skew version of Bollob\'as' Theorem for subspaces.

\begin{thm} \label{Lovasz}
Let $\F$ be an arbitrary field.  Let  $U_1, \ldots U_m$ be $r$-dimensional and $V_1, \ldots V_m$ be $s$-dimensional subspaces of a linear  space $W$ over the field $\F$. Assume that 
\begin{itemize}
\item[(i)] $U_i \cap V_i =\{0\}$ for each $1\leq i \leq m$;
\item[(ii)] $U_i\cap V_j\ne \{0\}$ whenever $i< j$ ($1\leq i, j \leq m$).
\end{itemize}
Then
$$
m\leq {r+s \choose r}.
$$
\end{thm}

In this paper our main aim is to give a subspace version of Theorem \ref{Tuza1} and \ref{Tuza2}.

The following definitions were motivated by Theorem \ref{Tuza2} and \ref{Lovasz}.

\begin{definition}
Let $\F$ be a fixed field. We say that a system $\cS=\{(U_i,V_i):~ 1\leq i\leq m\}$ is a  {\em weak} $ISP$-system of subspaces of an $n$-dimensional  linear space $W$ over the field ${\F}$, if $\cS$ satisfies the following  conditions:
\begin{itemize}
\item[(i)] $U_i \cap V_i =\{0\}$ for each $1\leq i \leq m$;
\item[(ii)] $U_i\cap V_j\ne \{0\}$ or $U_j\cap V_i\ne \{0\}$ for $i\ne j$ ($1\leq i, j \leq m$).
\end{itemize}
\end{definition}

\begin{definition}
Let $\F$ be a fixed field. We say that a  system $\cS=\{(U_i,V_i):~ 1\leq i\leq m\}$ of subspaces of a linear space $W$ over the field ${\F}$ is a {\em weak $(u,v)$-system}, if $\cS$ satisfies the conditions
\begin{itemize}
\item[(i)] $\cS$ is a   weak $ISP$-system;
\item[(ii)] $dim(U_i)=u$ \mbox{ and } $dim(V_i)=v$ for each $1\leq i \leq m$.
\end{itemize}
\end{definition}

Our main results are upper bounds for the size of weak $ISP$-systems and weak $(u,v)$-systems. 

\begin{thm} \label{main2}
Let $\cS=\{(U_i,V_i):~ 1\leq i\leq m\}$ be a weak $ISP$-system of subspaces of a linear space $W$ over the finite field ${\Fq}$. Let $u_i:=dim(U_i)$ and $v_i:=dim(V_i)$ for each $1\leq i\leq m$. Let $0\leq j\leq n$ be an arbitrary, but fixed integer. Then we have
$$
\sum_{i=1}^m \frac{\Big[ { n-v_i-u_i \atop j-u_i} \Big]_q  q^{(j-u_i)v_i}}{\Big[ { n \atop j} \Big]_q}\leq 1.
$$
\end{thm}

\begin{thm} \label{main3}
Let $\cS=\{(U_i,V_i):~ 1\leq i\leq m\}$ be a weak $(u,v)$-system of subspaces of an $n$-dimensional linear space $W$ over the finite field $\Fq$. Then
$$
m\leq \left(\frac{q}{q-1} \right)^n q^{uv}.
$$
\end{thm}

\section{Proofs of our main results}

In the proof of our main results we use  the following bounds for the $q$-binomial coefficients.
\begin{lemma} \label{qbinom}
Let $0\leq j\leq n$ be natural numbers. Then
$$
\Big[ { n \atop j} \Big]_q \leq \Big( \frac{q}{q-1} \Big)^n q^{j(n-j)}.
$$
\end{lemma}
\proof
This follows immediately from the inequalities 
$$
q^{{n\choose 2}}\leq [n]_q!\leq\Big( \frac{q}{q-1} \Big)^n q^{{n\choose 2}}.
$$
\qed

We use also in the proof of Theorem \ref{main2} the following simple Lemma (see Lemma 2.2  in \cite{WL}).

\begin{lemma} \label{disjoint}
Let $V$ denote the $n$-dimensional vector space over the finite field  ${\Fq}$ and fix an  $(n-d)$-dimensional subspace $K$ of $V$, where $0\leq d\leq n$. Let  $U_1$  be a fixed  $\ell_1$-subspace of $V$   such that $U_1\cap K =\{0\}$.  Let $u(n,d; \ell_1, \ell_2)$ denote 
the number of $\ell_2$-subspaces $U_2$ of $V$ satisfying $U_2 \cap K=\{0\}$ and $U_1 \subseteq U_2$. Then
$$
u(n,d; \ell_1, \ell_2)=\frac{\Big[ { d \atop \ell_2} \Big]_q\Big[ { \ell_2 \atop \ell_1} \Big]_q q^{(\ell_2-\ell_1)(n-d)}}{\Big[ { d \atop \ell_1} \Big]_q}.
$$ 
\end{lemma}
\qed

{\bf Proof of Theorem \ref{main2}:}

Let $1\leq i\leq m$, $0\leq  j\leq n$ be  fixed integers.   Let $\cF(i,j)$ denote the following subset of subspaces of $W$:
$$
\cF(i,j):= \{U\leq  W:~ \mbox{dim}(U)=j, U_i\subseteq U, V_i\cap U= \{0\}\}.
$$
Then it follows immediately from Lemma \ref{disjoint} that 
$$
|\cF(i,j)|= \frac{\Big[ { n-v_i \atop j} \Big]_q\Big[ { j \atop u_i} \Big]_q q^{(j-u_i)v_i}}{\Big[ { n-v_i \atop u_i} \Big]_q}.
$$
for each  $0\leq j\leq n$.

\begin{lemma} \label{disj}
Let $0\leq j\leq n$ be fixed. Let $1\leq i_1 < i_2 \leq m$ be two indices. Then 
$$
\cF(i_1,j)\cap \cF(i_2,j)=\emptyset.
$$
\end{lemma}

\proof
We can prove this statement by an indirect argument. Suppose that there exist two indices   $1\leq i_1 < i_2 \leq m$ such that $\cF(i_1,j)\cap \cF(i_2,j)\ne \emptyset$. Let $U\in \cF(i_1,j)\cap \cF(i_2,j)$ be an arbitrary, but fixed subspace. Then $U_{i_1}\subseteq U$ and $V_{i_1}\cap U=\{0\}$. Similarly 
 $U_{i_2}\subseteq U$ and $V_{i_2}\cap U=\{0\}$. Hence we get that
$$
U_{i_1} \cap V_{i_2}=\{0\}
$$
and 
$$
U_{i_2} \cap V_{i_1}=\{0\},
$$
which gives a contradiction, because $\cS=\{(U_i,V_i):~ 1\leq i\leq m\}$ is a weak $(u,v)$-system of subspaces of the linear space $W$. \qed

In the following let $0\leq j\leq n$ be a fixed integer.

It follows from Lemma \ref{disj} that 
$$                                                  
\sum_{i=1}^m |\cF(i,j)|=|\bigcup\limits_{i=1}^m \cF(i,j)|\leq \Big[ { n \atop j} \Big]_q,
$$
because $\cF(i,j)\subseteq \{U\leq W:~ \mbox{dim}(U)=j\}$.
Hence
\begin{equation}\label{qbinom2}
\sum_{i=1}^m \frac{\Big[ { n-v_i \atop j} \Big]_q\Big[ { j \atop u_i} \Big]_q q^{(j-u_i)v_i}}{\Big[ { n-v_i \atop u_i} \Big]_q}
\leq \Big[ { n \atop j} \Big]_q
\end{equation}
But it is easy to verify that
$$
\frac{\Big[ { n-v_i \atop j} \Big]_q\Big[ { j \atop u_i} \Big]_q }{\Big[ { n-v_i \atop u_i} \Big]_q}=\Big[ { n-v_i-u_i \atop j-u_i} \Big]_q,
$$
hence it follows from inequality (\ref{qbinom2}) that 
$$
\sum_{i=1}^m \Big[ { n-v_i-u_i \atop j-u_i} \Big]_q  q^{(j-u_i)v_i} \leq \Big[ { n \atop j} \Big]_q,
$$
which was to be proved. \qed

{\bf Proof of Theorem \ref{main3}:} 
If $\cS=\{(U_i,V_i):~ 1\leq i\leq m\}$ is a weak $(u,v)$-system of subspaces of the linear space $W$, then  $u_i=dim(U_i)=u$ and $v_i=dim(V_i)=v$ for each $1\leq i\leq m$. It follows from Theorem \ref{main2} that
$$
\sum_{i=1}^m \frac{\Big[ { n-u-v \atop j-u} \Big]_q  q^{(j-u)v}}{\Big[ { n \atop j} \Big]_q}\leq 1
$$
for each $1\leq j\leq n$. 
Let $j:=n-v$. This choice implies that
$$
\sum_{i=1}^m \frac{  q^{(n-v-u)v}}{\Big[ { n \atop v} \Big]_q}\leq 1.
$$
It follows from Lemma \ref{qbinom} that
$$               
\sum_{i=1}^m \frac{  q^{(n-v-u)v}}{ \Big( \frac{q}{q-1} \Big)^n q^{v(n-v)} }\leq 1.
$$
But then 
$$
 m\frac{  q^{-uv}}{ {\Big( {\frac{q}{q-1}} \Big)}^n }\leq 1,
$$
which was to be proved. \qed

\section{Concluding remarks}

We can raise  the following natural question: Let $u,v$ be fixed positive integers. Let $\F$ be a fixed field. Determine the largest integer $t:=t(u,v)$ such that
there exists a weak $(u,v)$-system  $\cS=\{(U_i,V_i):~ 1\leq i\leq t\}$  of $t(u,v)$ pairs 
of subspaces of an $n$-dimensional linear space $W$ over the  field $\F$ .
 
If $\F$ is the finite field $\Fq$, then we proved in Theorem \ref{main3} that 
$$
t(u,v)\leq \left(\frac{q}{q-1} \right)^n q^{uv}.
$$

On the other hand, it is easy to verify the lower bound $m(u,v)\leq t(u,v)$. Namely let $\{e_1, \ldots ,e_n\}$ denote a fixed basis of the  $n$-dimensional linear space $W$ over $\F$. By the definition of the number $m(u,v)$ there exists a system  $\cS=\{(A_i,B_i):~ 1\leq i\leq m(u,v)\}$  of $m(u,v)$ pairs of sets satisfying the  conditions:
\begin{itemize}
\item[(i)] $A_1, \ldots A_m$ are $u$-element sets and  $B_1, \ldots B_m$ are $v$-element sets;
\item[(ii)] $A_i \cap B_i =\emptyset$ for each $1\leq i \leq m$;
\item[(iii)] $A_i\cap B_j\ne \emptyset$ or $A_j\cap B_i\ne  \emptyset$ for $i\ne j$ ($1\leq i, j \leq m$).
\end{itemize}

Define the generated subspaces $U_i:=\langle \{e_k:~ k\in A_i\}\rangle$ and  $V_i:=\langle \{e_l:~ l\in B_i\}\rangle$ for each $1\leq i\leq m(u,v)$. 

Then it is easy to verify that the system  $\cS=\{(U_i,V_i):~ 1\leq i\leq m(u,v)\}$  of $m(u,v)$ pairs 
of subspaces  is a weak $(u,v)$-system.


\begin{thebibliography}{MM}

\bibitem{B} B. Bollob\'as,  On generalized graphs. {\em Acta Math. Hung.} {\bf 16(3)}, (1965) 447-452.

\bibitem{BF} L. Babai and P. Frankl, {\em Linear algebra methods in
combinatorics}, September 1992.

\bibitem{K} W. Koepf,  {\em Hypergeometric Summation: An Algorithmic Approach to Summation and Special Function Identities}, Vieweg, 1998.

\bibitem{L1} L. Lov\'asz, Flats in matroids and geometric graphs, in: {\em Combinatorial surveys}, Proc. 6th British Comb. Conf., Egham 1977, Acad. Press, London 1977, 45--86.

\bibitem{L2} L. Lov\'asz, Topological and algebraic methods in graph theory, in {\em Graph theory and related topics}, Proc. Conf., Univ. Waterloo, Waterloo, Ont., 1979,  1--14.

 \bibitem{T} Z. Tuza, Application of Set-Pair Method in Extremal Hypergraph Theory, in ``Extremal problems for Finite Sets'', {\em Bolyai Soc. Math. Studies} {\bf 3}, J\'anos Bolyai Math. Soc., Budapest, 1994, 479--514. 

\bibitem{T2} Z. Tuza,  Inequalities for two set systems with prescribed intersections. {\em Graphs and Comb.}, {\bf 3} (1),(1987) 75-80.

\bibitem{WL} K. Wang and Z. Li, Lattices associated with vector spaces over a finite field, {\em Lin. Algebra and Its Applications} {\bf 429(2)},  (2008), 439-446.

\end{thebibliography}
\end{document}